\documentclass[10pt]{amsart}

\usepackage{amssymb,amsmath,amsthm}

\usepackage{a4wide}

\usepackage{subfig}
\usepackage{hyphenat}
\usepackage[colorlinks=true]{hyperref}
\usepackage{mathtools}
\usepackage[numbers]{natbib}
\usepackage{enumerate}

\usepackage[capitalise]{cleveref}
\crefname{lemma}{Lemma}{Lemmas}
\crefname{equation}{Eq.}{Eqs.}
\Crefname{equation}{Equation}{Equations}
\crefname{assumption}{Assumption}{Assumptions}
\crefname{pluralequation}{Eqs.}{Eqs.}
\creflabelformat{pluralequation}{(#2#1#3)}

\creflabelformat{enumi}{#2#1#3)}

\numberwithin{equation}{section}

\newcommand{\1}{\mathbf{1}}
\newcommand{\C}{\mathbb{C}}
\newcommand{\EE}{\mathbb{E}}

\newcommand{\NN}{\mathbb{N}}
\newcommand{\Z}{\mathbb{Z}}
\newcommand{\R}{\mathbb{R}}
\newcommand{\HH}{\mathcal{H}}

\newcommand{\II}{\mathrm{II}}
\newcommand{\I}{\mathrm{I}}

\newcommand{\X}{\mathbf{X}}
\newcommand{\ZZ}{\mathbf{Z}}
\newcommand{\iii}{\mathrm{i}}
\newcommand{\dd}{\mathrm{d}}
\newcommand{\ee}{\mathrm{e}}
\newcommand{\ii}{\mathrm{i}}
\newcommand{\wt}{\widetilde}
\newcommand{\ub}{\underbrace}
\newcommand{\Leb}{\operatorname{Leb}}
\newcommand{\tr}{\operatorname{tr}}
\newcommand{\varop}{\operatorname{Var}}

\newcommand{\BSO}{\operatorname{B}}

\renewcommand{\leq}{\leqslant}

\theoremstyle{plain}
\newtheorem{theorem}{Theorem}[section]

\newtheorem{lemma}[theorem]{Lemma}

\theoremstyle{definition}

\theoremstyle{remark}
\newtheorem{remark}[theorem]{Remark}

\allowdisplaybreaks[1]

\usepackage{microtype}
\makeatletter
\def\MT@register@subst@font{\MT@exp@one@n\MT@in@clist\font@name\MT@font@list
   \ifMT@inlist@\else\xdef\MT@font@list{\MT@font@list\font@name,}\fi}
\makeatother 

\title{Eigenvalue distribution of large sample covariance matrices of linear processes}

\author{Oliver Pfaffel}
\author{Eckhard Schlemm}

\address{TUM Institute for Advanced Study \& Fakult\"at f\"ur Mathematik, Technische Universit\"at M\"unchen, Germany}

                                     %key words and phrases

               %2000 Mathematics Subject Classification

\begin{document}

\begin{abstract}
We derive the distribution of the eigenvalues of a large sample covariance matrix when the data is dependent in time. More precisely, the dependence for each variable $i=1,\ldots,p$ is modelled as a linear process $(X_{i,t})_{t=1,\ldots,n}=(\sum_{j=0}^\infty c_j Z_{i,t-j})_{t=1,\ldots,n}$, where $\{Z_{i,t}\}$ are assumed to be independent random variables with finite fourth moments. If the sample size $n$ and the number of variables $p=p_n$ both converge to infinity such that $y=\lim_{n\to\infty}{n/p_n}>0$, then the empirical spectral distribution of $p^{-1}\X\X^T$ converges to a non\hyp{}random distribution which only depends on $y$ and the spectral density of $(X_{1,t})_{t\in\Z}$. In particular, our results apply to (fractionally integrated) ARMA processes, which we illustrate by some examples.
\end{abstract}

\maketitle

\section{Introduction and main result}

A typical object of interest in many fields is the sample covariance matrix $(n-1)^{-1} \X\X^T$ of a data matrix $\X=(X_{i,t})_{it}$, $i=1,\ldots,p$, $t=1,\ldots,n$. The matrix $\X$ can be seen as a sample of size $n$ of $p$-dimensional data vectors. For fixed $p$ one can show, as $n$ tends to infinity, that under certain assumptions the eigenvalues of the sample covariance matrix converge to the eigenvalues of the true underlying covariance matrix \citep{Anderson2003}. However, the assumption $p\ll n$ may not be justified if one has to deal with high dimensional data sets, so that it is often more suitable to assume that the dimension $p$ is of the same order as the sample size $n$, that is
$p=p_n\to\infty$ such that 
\begin{align}
\lim_{n\to\infty}\frac{n}{p}\eqqcolon y\in(0,\infty).
\label{npy}
\end{align}
For a symmetric matrix $A$ with eigenvalues $\lambda_1,\ldots,\lambda_p$, we denote by
\begin{equation*}
F^A = \frac{1}{p} \sum_{i=1}^p \delta_{\lambda_i}
\end{equation*}
the spectral distribution of $A$, where $\delta_x$ denotes the Dirac measure located at $x$. This means that $pF^A(B)$ is equal to the number of eigenvalues of $A$ that lie in the set $B$. From now on we will call $p^{-1} \X\X^T$ the sample covariance matrix. Due to \cref{npy}, this change of normalization can be reversed by a simple transformation of the limiting spectral distribution. For notational convenience we suppress the explicit dependence of the occurring matrices on $n$ and $p$ where this does not cause ambiguity. 

The distribution of Gaussian sample covariance matrices of fixed size was first computed in \citep{wishart1928generalised}. Several years later, it was Marchenko and Pastur \citep{Marchenko1967} who considered the case where the random variables $\{X_{i,t}\}$ are more general i.\,i.\,d.\ random variables with finite second moments $\EE X_{11}^2=1$, and the number $p$ of variables is of the same order as the sample size $n$. They showed that the empirical spectral distribution (ESD) $F^{p^{-1} \X\X^T}$ of $p^{-1} \X\X^T$ converges, as $n\to\infty$, to a non\hyp{}random distribution $\hat{F}$, called limiting spectral distribution (LSD), given by
\begin{equation}
\label{eq-MPlaw}
\hat{F}(\dd x) = \frac{1}{2\pi x} \sqrt{(x_+-x)(x-x_{-})} \1_{\{x_-\leq x\leq x_+\}} \dd x,
\end{equation}
and point mass $\hat F(\{0\})=1-y$ if $y<1$; in this formula, $x_{\pm}=(1\pm\sqrt{y})^2$. Here and in the following, convergence of the ESD means almost sure convergence as a random element of the space of probability measures on $\R$ equipped with the weak topology. In particular, the eigenvalues of the sample covariance matrix of a matrix with independent entries do not converge to the eigenvalues of the true covariance matrix, which is the identity matrix and therefore only has eigenvalue one. This leads to the failure of statistics that rely on the eigenvalues of $p^{-1} \X\X^T$ which have been derived under the assumption of fixed $p$, and random matrix theory is a tool to correct these statistics \citep{bai2009corrections,Johnstone2001}. In the case where the true covariance matrix is not the identity matrix, the LSD can in general only be given in terms of a non\hyp{}linear equation for its Stieltjes transform, which is defined by
\begin{equation*}
m_{\hat{F}}(z) = \int \frac{1}{\lambda-z} \dd\hat{F} \quad \forall z\in\C^+\coloneqq\{z=u+\iii v\in\C: \Im z=v>0\}.
\end{equation*}
Conversely, the distribution $\hat{F}$ can be obtained from its Stieltjes transform $m_{\hat F}$ via the Stieltjes--Perron inversion formula (\citep[Theorem B.8]{bai2010}), which states that
\begin{equation}
\label{eq-Stieltjes}
\hat{F}([a,b])=\frac{1}{\pi}\lim_{\epsilon\to 0^+} \int_a^b \Im m_{\hat{F}}(x+ \iii\epsilon) \dd x.
\end{equation}
for all continuity points $a<b$ of $\hat F$. For a comprehensive account of random matrix theory we refer the reader to \citep{Anderson2009,bai2010,Mehta2004}, and the references therein. 

Our aim in this paper is to obtain a Marchenko--Pastur type result in the case where there is dependence within the rows of $\X$. More precisely, for $i=1,\ldots,p$, the $i$th row of $\X$ is given by a linear process of the form
\begin{equation*}
(X_{i,t})_{t=1,\ldots,n} = \left(\sum_{j=0}^\infty c_j Z_{i,t-j}\right)_{t=1,\ldots,n}, \quad c_j\in\R.
\end{equation*} 
Here,  $(Z_{i,t})_{it}$ is an array of independent random variables that satisfies
\begin{align}
\EE Z_{i,t}=0,\quad \EE Z_{i,t}^2=1, \quad\text{ and }\quad \sigma_4\coloneqq\sup_{i,t}\EE Z_{i,t}^4<\infty,
\label{Z1}
\end{align}
as well as the Lindeberg\hyp{}type condition that, for each $\epsilon>0$,
\begin{align}
\frac{1}{pn} \sum_{i=1}^p \sum_{j=1}^n \EE\left(Z_{i,t}^2 \1_{\{Z_{i,t}^2\geqslant \epsilon n\}}\right) \to 0,\quad \text{ as }\quad n\to\infty.
\label{Z2}
\end{align}
Clearly, \cref{Z2} is satisfied if all $\{Z_{i,t}\}$ are identically distributed.

The novelty of our result is that we allow for dependence within the rows, and that the equation for $m_{\hat{F}}$ is given in terms of the spectral density
\begin{equation*}
f(\omega) = \sum_{h\in\Z} \gamma(h) \ee^{-\ii h\omega}, \quad \omega\in[0,2\pi],
\end{equation*}
of the linear processes $X_i$ only, which is the Fourier transform of the autocovariance function
\begin{equation*}
\gamma(h) = \sum_{j=0}^\infty c_j c_{j+|h|}, \quad h\in\Z.
\end{equation*}
Potential applications arise whenever data is not independent in time such that the Marchenko--Pastur law is not a good approximation. This includes e.\,g.\ wireless communications \citep{Tulino2004} and mathematical finance \citep{Bouchaud2009,Plerou2002}. Note that a similar question is also discussed in \citep{bai2008}. However, they have a different proof which relies on a moment condition to be verified. Furthermore, they assume that the random variables $\{Z_{i,t}\}$ are identically distributed so that the processes within the rows are independent copies of each other. More importantly, their results do not yield concrete formulas except in the AR(1) case and are therefore not directly applicable. In the context of free probability theory, the limiting spectral distribution of large sample covariance matrices of Gaussian ARMA processes is investigated in \citep{burda2010random}.

Before we present the main result of this article, we explain the notation used in this article. The symbols $\Z$, $\NN$ $\R$, and $\C$ denote the sets of integers, natural, real, and complex numbers, respectively. For a matrix $A$, we write $A^T$ for its transpose and $\tr A$ for its trace. Finally, the indicator of an expression $\mathcal{E}$ is denoted by $I_{\{\mathcal{E}\}}$ and defined to be one if $\mathcal{E}$ is true, and zero otherwise; for a set $S$, we also write $I_S(x)$ instead of $I_{\{x\in S\}}$.

\begin{theorem}
\label{maintheorem}
For each $i=1,\ldots, p$, let $X_{i,t}=\sum_{j=0}^\infty{c_j Z_{i,t-j}}$, $t\in\Z$, be a linear stochastic process with continuously differentiable spectral density $f$. Assume that
\begin{enumerate}[i)]
  \item\label{maintheoremconditionsZ} the array $(Z_{i,t})_{it}$ satisfies conditions \labelcref{Z1,Z2},
  \item\label{maintheoremconditionsc} there exist positive constants $C$ and $\delta$ such that $|c_j|\leq C(j+1)^{-1-\delta}$ for all $j\geqslant 0$,
  \item\label{maintheoremconditionsflevelsets} for almost all $\lambda\in\R$, $f(\omega)=\lambda$ for at most finitely many $\omega\in[0,2\pi]$, and
  \item\label{maintheoremconditionsdfzero} $f'(\omega)\neq 0$ for almost every $\omega$.
\end{enumerate}
Then the empirical spectral distribution $F^{p^{-1}\X\X^T}$ of $p^{-1}\X\X^T$ converges, as $n$ tends to infinity, almost surely to a non\hyp{}random probability distribution $\hat{F}$ with bounded support. Moreover, there exist positive numbers $\lambda_-,\lambda_+$ such that the Stieltjes transform $z\mapsto m_{\hat F}(z)$ of $\hat F$ is the unique mapping $\C^+\to\C^+$ satisfying
%unique solution in $\C^+=\{z\in\C: \Im(z)>0\}$ of
\begin{equation}
\label{eq-stieltjes}
\frac{1}{m_{\hat F}(z)} = - z +  \frac{y}{2\pi} \int{\lambda_{-}}^{\lambda_{+}}{\frac{\lambda}{1+\lambda m_{\hat F}(z)} \sum_{\omega\in[0,2\pi]:f(\omega)=\lambda} \frac{1}{\left|f'(\omega)\right|} \dd\lambda}.
\end{equation}
\end{theorem}

The assumptions of the theorem are met, for instance, if $(X_{i,t})_t$ is an ARMA or fractionally integrated ARMA process; see \cref{examples} for details.

\Cref{maintheorem}, as it stands, does not contain the classical Marchenko--Pastur law as a special case. For if the entries $X_{i,t}$ of the matrix $\X$ are i.\,i.\,d., the corresponding spectral density $f$ is identically equal to the variance of $X_{1,1}$, and thus condition \labelcref{maintheoremconditionsdfzero} is not satisfied. We therefore also present a version of \cref{maintheorem} that holds if the rows of the matrix $\X$ have a piecewise constant spectral density.
\begin{theorem}
\label{maintheorempiecewise}
For each $i=1,\ldots, p$, let $X_{i,t}=\sum_{j=0}^\infty{c_j Z_{i,t-j}}$, $t\in\Z$, be a linear stochastic process with spectral density $f$ of the form
\begin{equation}
f:[0,2\pi]\to\R^+,\quad \omega \mapsto \sum_{j=1}^k{\alpha_j \1_{A_j}(\omega)},\quad k\in\NN,
\end{equation}
for some positive real numbers $\alpha_j$ and a measurable partition $A_1\cup\cdots\cup A_k$ of the interval $[0,2\pi]$. If conditions \labelcref{maintheoremconditionsZ,maintheoremconditionsc} of \cref{maintheorem} hold, then the empirical spectral distribution $F^{p^{-1}\X\X^T}$ of $p^{-1}\X\X^T$ converges, as $n\to\infty$, almost surely to a non\hyp{}random probability distribution $\hat{F}$ with bounded support. Moreover, the Stieltjes transform $z\mapsto m_{\hat F}(z)$ of $\hat F$ is the unique mapping $\C^+\to\C^+$ that satisfies
%unique solution in $\C^+=\{z\in\C: \Im(z)>0\}$ of
\begin{equation}
\label{eq-stieltjespiecewise}
\frac{1}{m_{\hat F}(z)} = - z +  \frac{y}{2\pi} \sum_{j=1}^k \frac{|A_j|\alpha_j}{1+\alpha_j m_{\hat F}(z)},
\end{equation}
where $|A_j|$ denotes the Lebesgue measure of the set $A_j$. In particular, if the entries of $\X$ are i.\,i.\,d.\ with unit variance, one recovers the limiting spectral distribution \labelcref{eq-MPlaw} of the Marchenko--Pastur law.
\end{theorem}

\begin{remark}
In applications one often considers processes of the form $X_{i,t} = \mu + \sum_{j=0}^\infty c_j Z_{i,t-j}$ with mean $\mu\neq 0$. If we denote by $x_t\in\R^p$ the $t$th column of the matrix $\X$, and define the empirical mean by  $\overline{x}=p^{-1}\sum_{t=1}^n x_{t}$, then the sample covariance matrix is given by the expression $p^{-1} \sum_{t=1}^n (x_t-\overline{x})(x_t-\overline{x})^T$ instead of $p^{-1}\X\X^T$. However, by \citep[Theorem A.44]{bai2010}, the subtraction of the empirical mean does not change the LSD, and thus \cref{maintheorem,maintheorempiecewise} remain valid if the underlying linear process has a non\hyp{}zero mean. 
%if the assumption $\EE Z_{i,t}=0$ is dropped.
\end{remark}
\begin{remark}
The proof of \cref{maintheorem,maintheorempiecewise} can easily be generalized to cover non\hyp{}causal linear processes, which are defined as $X_{i,t}=\sum_{j=-\infty}^\infty{c_j Z_{i,t-j}}$. For this case one obtains the same result except that the autocovariance function is now given by $\sum_{j=-\infty}^\infty c_j c_{j+|h|}$.
\end{remark}
\begin{remark}
If one considers a matrix $\X$ which has independent linear processes in its columns instead of its rows, one obtains the same formulas as in \cref{maintheorem,maintheorempiecewise} except that $y$ is replaced by $y^{-1}$. This is due to the fact that $\X^T\X$ and $\X\X^T$ have the same non\hyp{}trivial eigenvalues.
\end{remark}

In \cref{proofs} we proceed with the proofs of \cref{maintheorem,maintheorempiecewise}. Thereafter we present some interesting examples in \cref{examples}.

\section{Proofs}
\label{proofs}
In this section we present our proofs of \cref{maintheorem,maintheorempiecewise}. Dealing with infinite\hyp{}order moving average processes directly is dfficult, and we therefore first prove a variant of these theorems for the truncated processes $\wt{X}_{i,t} = \sum_{j=0}^n c_j Z_{i,t-j}$. We define the $p\times n$ matrix $\wt{\X}=(\wt{X}_{i,t})_{it}$, $i=1,\ldots,p$, $t=1,\ldots,n$.
\begin{theorem}
\label{maintheoremmod}
Under the assumptions of \cref{maintheorem} (\cref{maintheorempiecewise}), the empirical spectral distribution of the sample covariance matrix of the truncated process $\wt X$ converges, as $n$ tends to infinity, to a deterministic distribution with bounded support. Its Stieltjes transform is uniquely determined by \cref{eq-stieltjes} (\cref{eq-stieltjespiecewise}).
\end{theorem}
\begin{proof}
The proof starts from the observation that one can write $\wt\X=\ZZ H$, where $\R^{p\times 2n}\ni\ZZ=(Z_{i,t})_{it}$, $i=1,\ldots,p$, $t=1-n,\ldots,n$, and 
\begin{align}
\label{eq-DefH}
H = \left(\begin{array}{ccccccccc}
      c_n 	& c_{n-1} 	& \ldots 	& c_1 	& c_0 	& 0	& \ldots& 0\\
      0		& c_n		& \ldots	& c_2	& c_1	& c_0	&	& \vdots\\
      \vdots	&		& \ddots	& \vdots& \vdots& 	& \ddots& 0\\
      0		& \ldots	& 0		& c_n	& c_{n-1}&\ldots& \ldots& c_0	
      \end{array}
\right)^T\in \R^{2n\times n}.
\end{align}
In particular, $\wt\X\wt\X^T=\ZZ HH^T\ZZ^T$. In order to prove convergence of the empirical spectral distribution $F^{p^{-1}\wt\X\wt\X^T}$ and to obtain a characterization of the limiting distribution, it suffices, by \citep[Theorem 1]{pan2010}, to prove that the spectral distribution $F^{HH^T}$ of $HH^T$ converges to a non\hyp{}trivial limiting distribution. This will be done in \cref{lemma-LSDHH}, where the LSD of $HH^T$ is shown to be $\hat F^{HH^T}=\frac{1}{2}\delta_0+\frac{1}{2}\hat F^\Gamma$; the distribution $\hat F^\Gamma$ is computed in \cref{lemma-densityGamma} if we impose the assumptions of \cref{maintheorem}, respectively in \cref{lemma-densityGammapiecewise} if we impose the assumptions of \cref{maintheorempiecewise}. Inserting this expression for $\hat F^{HH^T}$ into equation (1.2) of \citep{pan2010} shows that the ESD $F^{p^{-1}\wt\X\wt\X^T}$ converges, as $n\to\infty$, almost surely to a deterministic distribution, which is determined by the requirement that its Stieltjes transform $z\mapsto m(z)$ satisfies
\begin{equation}
\label{eq-stieltjescharacGamma}
 \frac{1}{m(z)} = - z +  2y \int_{\lambda_{-}}^{\lambda_{+}} \frac{\lambda }{1+\lambda m(z)}\dd\hat F^{HH^T}=- z +  y \int_{\lambda_{-}}^{\lambda_{+}} \frac{\lambda }{1+\lambda m(z)}\dd\hat F^\Gamma.
\end{equation}
Using the explicit formulas of $\hat F^\Gamma$ computed in \cref{lemma-densityGamma,lemma-densityGammapiecewise}, one obtains \cref{eq-stieltjes,eq-stieltjespiecewise}. Uniqueness of a mapping $m:\C^+\to\C^+$ solving \cref{eq-stieltjescharacGamma} was shown in \citep[p. 88]{bai2010}. We complete the proof by arguing that the LSD of $p^{-1}\wt\X\wt\X^T$ has bounded support. For this it is enough, by \citep[Theorem 6.3]{bai2010}, to show that the spectral norm of $HH^T$ is bounded in $n$, which is also done in \cref{lemma-LSDHH}.
\end{proof}

\begin{lemma}
 \label{lemma-LSDHH}
Let $H=(c_{n-i+j}\1_{\{0\leq n-i+j\leq n\}})_{ij}$ be the matrix appearing in \cref{eq-DefH}, and assume that there exist positive constants $C,\delta$ such that $|c_j|\leq C(j+1)^{-1-\delta}$ (assumption \labelcref{maintheoremconditionsc} of \cref{maintheorem}). Then the spectral norm of the matrix $HH^T$ is bounded in $n$. If, moreover, the spectral distribution of the Toeplitz matrix $\Gamma=(\gamma(i-j))_{ij}$ converges weakly to some limiting distribution $\hat F^{\Gamma}$, then the spectral distribution $F^{HH^T}$ converges weakly, as $n\to\infty$, to $\frac{1}{2}\delta_0+\frac{1}{2}\hat F^{\Gamma}$.
\end{lemma}

\begin{proof}
We first introduce the notation $\HH\coloneqq HH^T\in\R^{2n\times 2n}$ as well as the block decomposition $\HH=\left[\begin{array}{cc}\HH_{11} & \HH_{12}\\\HH_{12}^T & \HH_{22}\end{array}\right]$, $\HH_{ij}\in \R^{n\times n}$. We prove the second part of the lemma first. There are several ways to show that the spectral distributions of two sequences of matrices converge to the same limit. In our case it is convenient to use \citep[Corollary A.41]{bai2010} which states that two sequences $A_n$ and $B_n$, either of whose empirical spectral distribution converges, have the same limiting spectral distribution if $n^{-1}\tr(A_n-B_n)(A_n-B_n)^T$ converges to zero as $n$ tends to infinity. We shall employ this result twice: first to show that the LSDs of $\HH=HH^T$ and $\wt\HH\coloneqq\operatorname{diag}(0,\HH_{22})$ agree, and then to prove equality of the LSDs of $\HH_{22}$ and $\Gamma$. Let $\Delta_{\HH}=n^{-1}\tr(\HH-\widetilde\HH)(\HH-\widetilde\HH)^T$; a direct calculation shows that $\label{eq-DeltaHH}
\Delta_{\HH} = n^{-1}\left[\tr \HH_{11}\HH_{11}^T+2\tr \HH_{12}\HH_{12}^T\right]$, and we will consider each of the two terms in turn. From the definition of $H$ it follows that the $(i,j)$th entry of $\HH$ is given by $\HH^{ij}=\sum_{k=1}^n{c_{n-i+k}c_{n-j+k}\1_{\{\max{(i,j)}-n\leq k\leq \min{(i,j)}\}}}$. The trace of the square of the upper left block of $\HH$ therefore satisfies
\begin{align*}
\tr\HH_{11}\HH_{11}^T=\sum_{i,j=1}^n{\left\{\HH^{ij}\right\}^2} = & \sum_{i,j=1}^n\left[\sum_{k=1}^{\min{(i,j)}}{c_{n-i+k}c_{n-j+k}}\right]^2\\
  \leq & \sum_{i,j,k,l=1}^n{|c_{i+k-1}||c_{j+k-1}||c_{i+l-1}||c_{j+l-1}|}\\
  \leq & C^4 \sum_{i,j,k,l=2}^{n+1}{i^{-1-\delta}j^{-1-\delta}l^{-1-\delta}k^{-1-\delta}}\\
    < & \left[C\zeta(1+\delta)\right]^4<\infty,
\end{align*}
where $\zeta(z)$ denotes the Riemann zeta function. As a consequence, the limit of  $n^{-1}\tr\HH_{11}\HH_{11}^T$ as $n$ tends to infinity is zero. Similarly, we obtain for the trace of the square of the off-diagonal block of $\HH$ the bound
\begin{align*}
\tr\HH_{12}\HH_{12}^T=\sum_{i=1}^n{\sum_{j=n+1}^{2n}\left\{\HH^{ij}\right\}^2} =& \sum_{i=1}^n\sum_{j=n+1}^{n+i}\left[\sum_{k=j-n}^i{c_{n-i+k}c_{n-j+k}}\right]^2\\
  \leq & \sum_{i=1}^n\sum_{j=1}^{n}\sum_{k=j}^{n-i+1}\sum_{l=j}^{n-i+1}{c_{i+k-1}c_{k-j}c_{i+l-1}c_{l-j}}\\
  \leq & \sum_{i=1}^n\sum_{j=1}^n\sum_{r=0}^n\sum_{s=0}^n{|c_{i+r+j-1}||c_r||c_{s+j-1}||c_s|}\\
  \leq & C^4 \sum_{i,j,r,s=1}^{n+1}{ i^{-1-\delta} r^{-1-\delta} s^{-1-\delta} j^{-1-\delta}}\\
  < & \left[C\zeta(1+\delta)\right]^4<\infty,
\end{align*}
which shows that the limit of $n^{-1}\tr\HH_{12}\HH_{12}^T$ is zero. It follows that $\Delta_{\HH}$, as defined in \cref{eq-DeltaHH}, converges to zero as $n$ goes to infinity, and therefore that the LSDs of $\HH$ and  $\wt \HH=\operatorname{diag}(0,\HH_{22})$ coincide. The latter distribution is clearly given by $F^{\wt\HH}=\frac{1}{2}\delta_0+\frac{1}{2}F^{\HH_{22}}$, and we show next that the LSD of $\HH_{22}$ agrees with the LSD of $\Gamma=(\gamma(i-j))_{ij}$. As before it suffices to show, by \citep[Corollary A.41]{bai2010}, that $\Delta_{\Gamma}=n^{-1}\tr(\HH_{22}-\Gamma)(\HH_{22}-\Gamma)^T$ converges to zero as $n$ tends to infinity. It follows from the definitions of $\HH$ and $\Gamma$ that $n\Delta_{\Gamma}$ can be estimated as
\begin{align*}
n\Delta_{\Gamma} =& \sum_{i,j=1}^n\left[\sum_{k=\max{(i,j)}}^n{c_{k-i}c_{k-j}}-\sum_{k=1}^\infty{c_{k-1}c_{k+|i-j|-1}}\right]^2\\
  =&  \sum_{i,j=1}^n\left[\sum_{k=\max{(i,j)}}^n{c_{k-i}c_{k-j}}-\sum_{k=\max{(i,j)}}^\infty{c_{k-i}c_{k-j}}\right]^2\\
  =& \sum_{i,j=1}^n\sum_{k,l=1}^\infty{c_{k+i-1}c_{k+j-1}c_{l+i-1}c_{l+j-1}}\\
  \leq & C^4\sum_{i,j=2}^{n+1}\sum_{k,l=2}^\infty{i^{-1-\delta}j^{-1-\delta}k^{-1-\delta}l^{-1-\delta}}< \left[C\zeta(1+\delta)\right]^4<\infty.
\end{align*}
Consequently, $\Delta_{\Gamma}$ converges to zero as $n$ goes to infinity, and it follows that $\hat F^{\HH}=\frac{1}{2}\delta_0+\frac{1}{2}\hat F^{\Gamma}$.

In order to show that the spectral norm of $\HH=HH^T$ is bounded in $n$, we use Gerschgorin's circle theorem (\citep[Theorem 2]{gerschgorin1931abgrenzung}), which states that every eigenvalue of $\HH$ lies in at least one of the balls $B(\HH^{ii},R_i)$ with centre $\HH^i$ and radius $R_i$, $i=1,\ldots,2n$, where the radii $R_i$ are defined as $R_i = \sum_{j\neq i}\left|\HH^{ij}\right|$. We first note that the centres $\HH^{ii}$ satisfy
\begin{equation*}
\HH^{ii} = \sum_{k=\max\{1,i-n\}}^{\min\{i,n\}}c_{n-i+k}^2\leq\sum_{k=0}^n c_k^2 \leq \left[C\zeta(2+2\delta)\right]^2<\infty.
\end{equation*}
To obtain a uniform bound for the radii $R_i$ we first assume that $i=1,\ldots,n$. Then
\begin{align*}
\left|R_i\right| \leq& \sum_{j=1}^{n}\sum_{k=1}^{\min\{i,j\}}{|c_{n-i+k}||c_{n-j+k}|} + \sum_{j=n+1}^{2n}\sum_{k=j-n}^{i}{|c_{n-i+k}||c_{n-j+k}|}\\
  \leq & \sum_{j,k=1}^n{|c_{n-i+k}||c_{j+k-1}|} + \sum_{j=n+1-i}^{2n-i}\sum_{k=0}^{n-j}|c_{k+j}||c_k|\leq2\left[C\zeta(1+\delta)\right]^2<\infty.
\end{align*}
Similarly we find that, for $i=n+1,\ldots,2n$,
\begin{align*}
\left|R_i\right| \leq& \sum_{j=1}^n\sum_{k=i-n}^{j}{|c_{n-i+k}||c_{n-j+k}|} + \sum_{j=n+1}^{2n}\sum_{k=\max\{i,j\}-n}^{n}{|c_{n-i+k}||c_{n-j+k}|}\\
  \leq & \sum_{j=i-n}^{i-1}\sum_{k=0}^{n+1-j}|c_{k+j}||c_k| + \sum_{j=n+1}^{2n}\sum_{k=0}^{n-\max\{i,j\}}|c_{k}||c_{k+|j-i|}|\leq3\left[C\zeta(1+\delta)\right]^2
\end{align*}
is bounded, which completes the proof.
\end{proof}

In the following two lemmas, we argue that the distribution $\hat F^\Gamma$ exists and we prove explicit formulas for it in the case that the assumptions of \cref{maintheorem} or \cref{maintheorempiecewise} are satisfied.

\begin{lemma}
\label{lemma-densityGamma}
Let $(c_j)_j$ be a sequence of real numbers, $\gamma:h\mapsto\sum_{j=0}^\infty{c_jc_{j+|h|}}$, and $f:\omega\mapsto\sum_{h\in \Z}{\gamma(h)}\ee^{-\ii h\omega}$. Under the assumptions of \cref{maintheorem} it holds that the spectral distribution $F^{\Gamma}$ of $\Gamma=(\gamma(i-j))_{ij}$ converges weakly, as $n\to\infty$, to an absolutely continuous distribution $\hat F^\Gamma$ with bounded support and density
\begin{equation}
\label{eq-densityGamma}
g:(\lambda_-,\lambda_+)\to\R^+,\quad \lambda\mapsto\frac{1}{2\pi}\sum_{\omega:f(\omega)=\lambda}\frac{1}{\left|f'(\omega)\right|}.                 \end{equation}
\end{lemma}

\begin{proof}
We first note that under assumption \labelcref{maintheoremconditionsc} of \cref{maintheorem} the autocovariance function $\gamma$ is absolutely summable because
\begin{equation*}
\sum_{h=0}^\infty|\gamma(h)| \leq \sum_{h=0}^\infty\sum_{j=0}^\infty{|c_j||c_{j+h}|}\leq C^2\sum_{h,j=1}^\infty{h^{-1-\delta}j^{-1-\delta}}<\left[C\zeta(1+\delta\right]^2<\infty.
\end{equation*}
Szeg{\H{o}}'s first convergence theorem (\citep{grenander1984} and \citep[Corollary 4.1]{gray2006}) then implies that $\hat F^\Gamma$ exists, and that the cumulative distribution function of the eigenvalues of the Toeplitz matrix $\Gamma$ associated with the sequence $h\mapsto\gamma(h)$ is given by
\begin{equation}
\label{eq-CDFGamma}
G(\lambda)\coloneqq\frac{1}{2\pi}\int_0^{2\pi}\1_{\{f(\omega)\leq\lambda\}}\dd \omega=\frac{1}{2\pi}\Leb(\{\omega\in[0,2\pi]:f(\omega)\leq\lambda\}),
\end{equation}
for all $\lambda$ such that the level sets $\{\omega\in[0,2\pi]:f(\omega)=\lambda\}$ have Lebesgue measure zero. By assumption \labelcref{maintheoremconditionsflevelsets} of \cref{maintheorem}, \cref{eq-CDFGamma} holds for almost all $\lambda$. In order to prove that the LSD $\hat F^{\Gamma}$ is absolutely continuous with respect to the Lebesgue measure, it suffices to prove that the cumulative distribution function $G$ is differentiable almost everywhere. Clearly, for $\Delta\lambda>0$,
\begin{equation*}
G(\lambda+\Delta\lambda)-G(\lambda) = \frac{1}{2\pi}\Leb(\{\omega\in[0,2\pi]:\lambda<f(\omega)\leq\lambda+\Delta\lambda\}). 
\end{equation*}
Due to assumption \labelcref{maintheoremconditionsdfzero} of \cref{maintheorem}, the set of all $\lambda\in\R$ such that the set $\{\omega:\in[0,2\pi]:f(\omega)=\lambda \textnormal{ and } f'(\omega)=0\}$ is non\hyp{}empty is a Lebesgue null-set. Hence it is enough to consider only $\lambda$ for which this set is empty.
Let $f^{-1}(\lambda)=\{\omega: f(\omega)=\lambda\}$ be the pre\hyp{}image of $\lambda$, which is a finite set by assumption \labelcref{maintheoremconditionsflevelsets}. The implicit function theorem then asserts that, for every $\omega\in f^{-1}(\lambda)$, there exists an open interval $I_\omega$ around $\omega$ such that $f$ restricted to $I_\omega$ is invertible. It is no restriction to assume that these $I_\omega$ are disjoint. By choosing $\Delta\lambda$ sufficiently small it can be ensured that the interval $[\lambda,\Delta\lambda]$ is contained in $\bigcap_{\omega\in f^{-1}(\lambda)}f(I_{\omega})$, and from the continuity of $f$ it follows that outside of $\bigcup_{\omega\in f^{-1}(\lambda)}I_\omega$, the values of $f$ are bounded away from $\lambda$, so that
\begin{align*}
&\lim_{\Delta\lambda\to 0}\frac{1}{\Delta\lambda}\left[G(\lambda+\Delta\lambda)-G(\lambda)\right]\\
=& \frac{1}{2\pi}\lim_{\Delta\lambda\to 0}\frac{1}{\Delta\lambda}\Leb\left(\bigcup_{\omega\in f^{-1}(\lambda)}\{\omega'\in I_\omega: \lambda<f(\omega')\leq\lambda+\Delta\lambda\}\right)\\
  =& \frac{1}{2\pi}\sum_{\omega\in f^{-1}(\lambda)}\lim_{\Delta\lambda\to 0}\frac{1}{\Delta\lambda}\Leb\left(\{\omega'\in I_\omega: \lambda<f(\omega')\leq\lambda+\Delta\lambda\}\right).
\end{align*}
In order to further simplify this expression, we denote the local inverse functions by $f^{-1}_\omega: f(I_\omega)\to[0,2\pi]$. Observing that the Lebesgue measure of an interval is given by its length, and that the derivatives of $f_\omega^{-1}$ are given by the inverse of the derivative of $f$, it follows that
\begin{align*}
\lim_{\Delta\lambda\to 0}\frac{1}{\Delta\lambda}\left[G(\lambda+\Delta\lambda)-G(\lambda)\right] =&\frac{1}{2\pi}\sum_{\omega\in f^{-1}(\lambda)}\lim_{\Delta\lambda\to 0}\frac{1}{\Delta\lambda}\left|f_\omega^{-1}(\lambda+\Delta\lambda)-f_\omega^{-1}(\lambda)\right|\\
  =& \frac{1}{2\pi}\sum_{\omega\in f^{-1}(\lambda)}\left|\frac{\dd}{\dd\lambda}f_\omega^{-1}(\lambda)\right|\\
  =& \frac{1}{2\pi}\sum_{\omega\in f^{-1}(\lambda)}\frac{1}{\left|f'(\omega)\right|}.
\end{align*}
This shows that $G$ is differentiable almost everywhere with derivative $g:\lambda\mapsto \frac{1}{2\pi}\sum_{\omega\in f^{-1}(\lambda)}\frac{1}{\left|f'(\omega)\right|}$. It remains to argue that the support of $\hat F^\Gamma$ is bounded. The absolute summability of $\gamma(\cdot)$ implies boundedness of its Fourier transform $f$. The claim then follows from \cref{eq-CDFGamma}, which shows that the support of $g$ is equal to the range of $f$.
\end{proof}

\begin{lemma}
\label{lemma-densityGammapiecewise}
Let $f:\omega\mapsto \sum_{j=1}^k \alpha_j \1_{A_j}(\omega)$ be the piecewise constant spectral density of the linear process $X_t=\sum_{j=0}^\infty c_j Z_{t-j}$, and denote the corresponding autocovariance function by $\gamma:h\mapsto\sum_{j=0}^\infty c_j c_{j+|h|}$. Under the assumptions of \cref{maintheorempiecewise} it holds that the spectral distribution $F^{\Gamma}$ of $\Gamma=(\gamma(i-j))_{ij}$ converges weakly, as $n\to\infty$, to the distribution $\hat F^\Gamma=(2\pi)^{-1}\sum_{j=1}^k{|A_j|\delta_{\alpha_j}}$.
\end{lemma}
\begin{proof}
Without loss of generality we may assume that $0<\alpha_1<\ldots<\alpha_k$. As in the proof of \cref{lemma-densityGamma} one sees that $\hat F^\Gamma$ exists, and that $\hat F^\Gamma(-\infty,\lambda)$ is given by
\begin{equation*}
\label{eq-CDFGammapiecewise}
G(\lambda)\coloneqq\frac{1}{2\pi}\Leb(\{\omega\in[0,2\pi]:f(\omega)\leq\lambda\}),\quad\forall \lambda\in[0,2\pi]\backslash \bigcup_{j=1}^k \{\alpha_j\}.
\end{equation*}
The special structure of $f$ thus implies that $G(\lambda)=(2\pi)^{-1}\sum_{j=1}^{k_\lambda}|A_j|$, where $k_\lambda$ is the largest integer such that $\alpha_{k_\lambda}\leq\lambda$. Since $G$ must be right\hyp{}continuous, this formula holds for all $\lambda$ in the interval $[0,2\pi]$. It is easy to see that the function $G$ is the cumulative distribution function of the discrete measure $(2\pi)^{-1}\sum_{j=1}^k{|A_j|\delta_{\alpha_j}}$, which completes the proof.
\end{proof}

\begin{proof}{of Theorems \mbox{\labelcref{maintheorem}} and \mbox{\labelcref{maintheorempiecewise}}}
It is only left to show that the truncation performed in \cref{maintheoremmod} does not alter the LSD, i.\,e.\ that the difference of $F^{p^{-1}\X\X^T}$ and $F^{p^{-1}\wt{\X}\wt{\X}^T}$ converges to zero almost surely. By \citep[Corollary A.42]{bai2010}, this means that we have to show that
\begin{equation}
\label{eq-IandII}
\ub{\frac{1}{p^2}\tr(\X\X^T+\wt{\X}\wt{\X}^T)}_{=\I} \ub{\frac{1}{p^2}\tr((\X-\wt{\X})(\X-\wt{\X})^T)}_{=\II}
\end{equation}
converges to zero. To this end we show that $\I$ has a limit, and that $\II$ converges to zero, both almost surely. By the definition of $\X$ and $\wt{\X}$ we have 
\begin{equation*}
\II = \frac{1}{p^2} \sum_{i=1}^p \sum_{t=1}^n \sum_{k=n+1}^\infty \sum_{m=n+1}^\infty c_k c_m Z_{i,t-k} Z_{i,t-m}.
\end{equation*}
We shall prove that the variances of $\II$ are summable. For this purpose we need the following two estimates which are implied by the Cauchy--Schwarz inequality, the assumption that $\sigma_4=\sup_{i,t}\EE Z_{i,t}^4$ is finite, and the assumed absolute summability of the coefficients $(c_j)_j$:
\begin{subequations}
\label{eq-fubini}
\begin{equation}
\label{eq-fubini1}
\EE\sum_{i=1}^p\sum_{t=1}^n\sum_{k,m=1}^\infty|c_k c_m Z_{i,t-k} Z_{i,t-m}|\leq pn\left(\sum_{k=1}^\infty |c_k|\right)^2 <\infty,
\end{equation}
\begin{align}
\label{eq-fubini2}
&\EE\sum_{i,i'=1=1}^p\sum_{t,t'=1}^n \sum_{k,k',m,m'=1}^\infty |c_k c_m c_{k'} c_{m'}  Z_{i,t-k}  Z_{i,t-m} Z_{i',t'-k'} Z_{i',t'-m'}|\\
  \leq & (np)^2\sigma_4\left(\sum_{k=1}^\infty |c_k|\right)^4 <\infty.\notag
\end{align}
\end{subequations}
Therefore we can, by Fubini's theorem, interchange expectation and summation to bound the variance of $\II$ as
\begin{equation*}
\varop(\II) \leq \frac{1}{p^4} \sum_{i,i'=1}^p \sum_{t,t'=1}^n \sum_{\substack{k,k'\\m,m'}=n+1}^\infty c_k c_m c_{k'} c_{m'} \EE( Z_{i,t-k}  Z_{i,t-m} Z_{i',t'-k'} Z_{i',t'-m'}).
\end{equation*}
Considering separately the terms where $i=i'$ and $i\neq i'$, we can write
\begin{align*}
\varop(\II)\leq&\frac{1}{p^4}\sum_{\substack{i,i'=1\\i\neq i'}}^p \sum_{t,t'=1}^n \sum_{\substack{k,k'\\m,m'}=n+1}^\infty c_k c_m c_{k'} c_{m'} \EE( Z_{i,t-k}  Z_{i,t-m} Z_{i',t'-k'} Z_{i',t'-m'})\\
  & + \frac{1}{p^4}\sum_{i=1}^p \sum_{t,t'=1}^n\sum_{\substack{k,k'\\m,m'}=n+1}^\infty c_k c_m c_{k'} c_{m'} \EE( Z_{i,t-k}  Z_{i,t-m} Z_{i,t'-k'} Z_{i,t'-m'}).
\end{align*}
For the expectation in the first sum not to be zero, $k$ must equal $m$ and $k'$ must equal $m'$, in which case its value is unity. The expectation in the second term can always be bounded by $\sigma_4$, so that we obtain
\begin{equation*}
\varop(\II) \leq \frac{p^2-p}{p^4}n^2\left(\sum_{k=n+1}^\infty c_k^2\right)^2+\sigma_4\frac{pn^2}{p^4}\left(\sum_{k=n+1}^\infty |c_k|\right)^4.
\end{equation*}
Due to \cref{npy} and the assumed polynomial decay of $c_k$ there exists a constant $K$ such that the right hand side is bounded by $Kn^{-1-4\delta}$, which implies that
\begin{equation*}
\sum_{n=1}^\infty \varop{(\II)} \leq K\sum_{n=1}^\infty n^{-1-4\delta}<\infty,
\end{equation*}
and therefore, by the first Borel--Cantelli lemma, that $\II$ converges to a constant almost surely. In order to show that this constant is zero, it suffices to shows that the expectation of $\II$ converges to zero. Since $\EE Z_{i,t}=0$, and the $\{Z_{i,t}\}$ are independent, one sees, using \cref{eq-fubini1} and again Fubini's theorem, that $\EE(\II)=n p^{-1} \sum_{k=n+1}^\infty c_k^2$, which converges to zero because the $\{c_k\}$ are square\hyp{}summable.

We now consider factor $\I$ of expression \labelcref{eq-IandII} and define $\Delta_X=\X\X^T-\wt{\X}\wt{\X}^T$. Then
\begin{equation}
\label{eq-IaandIb}
\I = \ub{\frac{1}{p^2}\tr(\Delta_X)}_{=\mathrm{I_a}} + 2 \ub{\frac{1}{p^2}\tr(\wt{\X}\wt{\X}^T)}_{=\mathrm{I_b}}.
\end{equation}
Because of
\begin{equation*}
(\X\X^T)_{ii} = \sum_{t=1}^n X_{i,t}^2 = \sum_{t=1}^n \sum_{k=0}^\infty \sum_{m=0}^\infty c_k c_m Z_{i,t-k} Z_{i,t-m},
\end{equation*}
and similarly $(\wt{\X}\wt{\X}^T)_{ii} = \sum_{t=1}^n \sum_{k=0}^n \sum_{m=0}^n c_k c_m Z_{i,t-k} Z_{i,t-m}$, we have that
\begin{align}
\tr(\Delta_X) =& \sum_{i=1}^p\left[ (\X\X^T)_{ii} - (\wt{\X}\wt{\X}^T)_{ii}\right]\notag\\
 = &\ub{\sum_{i=1}^p \sum_{t=1}^n \sum_{k=n+1}^\infty \sum_{m=n+1}^\infty c_k c_m Z_{i,t-k} Z_{i,t-m}}_{=\II\to 0\text{ a.s.}}\notag\\
\label{eq-trDelta}  & + 2\sum_{i=1}^p \sum_{t=1}^n \sum_{k=n+1}^\infty \sum_{m=1}^n c_k c_m Z_{i,t-k} Z_{i,t-m}.
\end{align}
\Cref{eq-fubini2} allows us to apply Fubini's theorem to compute the variance of the second term in the previous display as
\begin{equation*}
\frac{4}{p^4} \sum_{i,i'=1}^p \sum_{t,t'=1}^n \sum_{k,k'=n+1}^\infty \sum_{m,m'=1}^n c_k c_m c_{k'} c_{m'} \EE( Z_{i,t-k}  Z_{i,t-m} Z_{i',t'-k'} Z_{i',t'-m'}),
\end{equation*}
which is, by the same reasoning as we did for $\II$, bounded by
\begin{equation*}
4\sigma_4\frac{p}{p^4}n^2\left(\sum_{k=n+1}^\infty |c_k|\right)^2\left(\sum_{m=1}^n |c_m|\right)^2\leq Kn^{-1-2\delta},
\end{equation*}
for some positive constant $K$. Clearly, this is summable in $n$. Having, by \cref{eq-fubini1}, expected value zero, the second term of \cref{eq-trDelta} and, therefore, also $\tr(\Delta_X)$ both converge to zero almost surely. Thus, we only have to look at the contribution of $\mathrm{I_b}$ in expression \labelcref{eq-IaandIb}. From \cref{maintheoremmod} we know that $F^{p^{-1}\wt{\X}\wt{\X}^T}$ converges almost surely weakly to some non\hyp{}random distribution $\hat F$ with bounded support. Hence, denoting by $\lambda_1,\ldots,\lambda_p$ the eigenvalues of $p^{-1}\wt{\X}\wt{\X}^T$,
\begin{align*}
\mathrm{I_b} = \frac{1}{p} \tr\left(\frac{1}{p}\wt{\X}\wt{\X}^T\right) = \frac{1}{p} \sum_{i=1}^p \lambda_i 
= \int \lambda \dd F^{\frac{1}{p}\wt{\X}\wt{\X}^T} {\rightarrow} \int \lambda \dd\hat F < \infty,
\end{align*}
almost surely. It follows that, in \cref{eq-IandII}, factor $\I$ is bounded, and factor $\II$ converges to zero, and so the proof of \cref{maintheorem,maintheorempiecewise} is complete.
\end{proof}

\section{Illustrative examples}
\label{examples}
For several classes of widely employed linear processes, \cref{maintheorem} can be used to obtain an explicit description of the limiting spectral distribution. In this section we consider the class of autoregressive moving average (ARMA) processes as well as fractionally integrated ARMA models. The distributions we obtain in the case of AR(1) and MA(1) processes can be interpreted as one\hyp{}parameter deformations of the classical Marchenko--Pastur law.

\begin{figure}
\subfloat[$y=1$]{
\includegraphics[width=.31\linewidth]{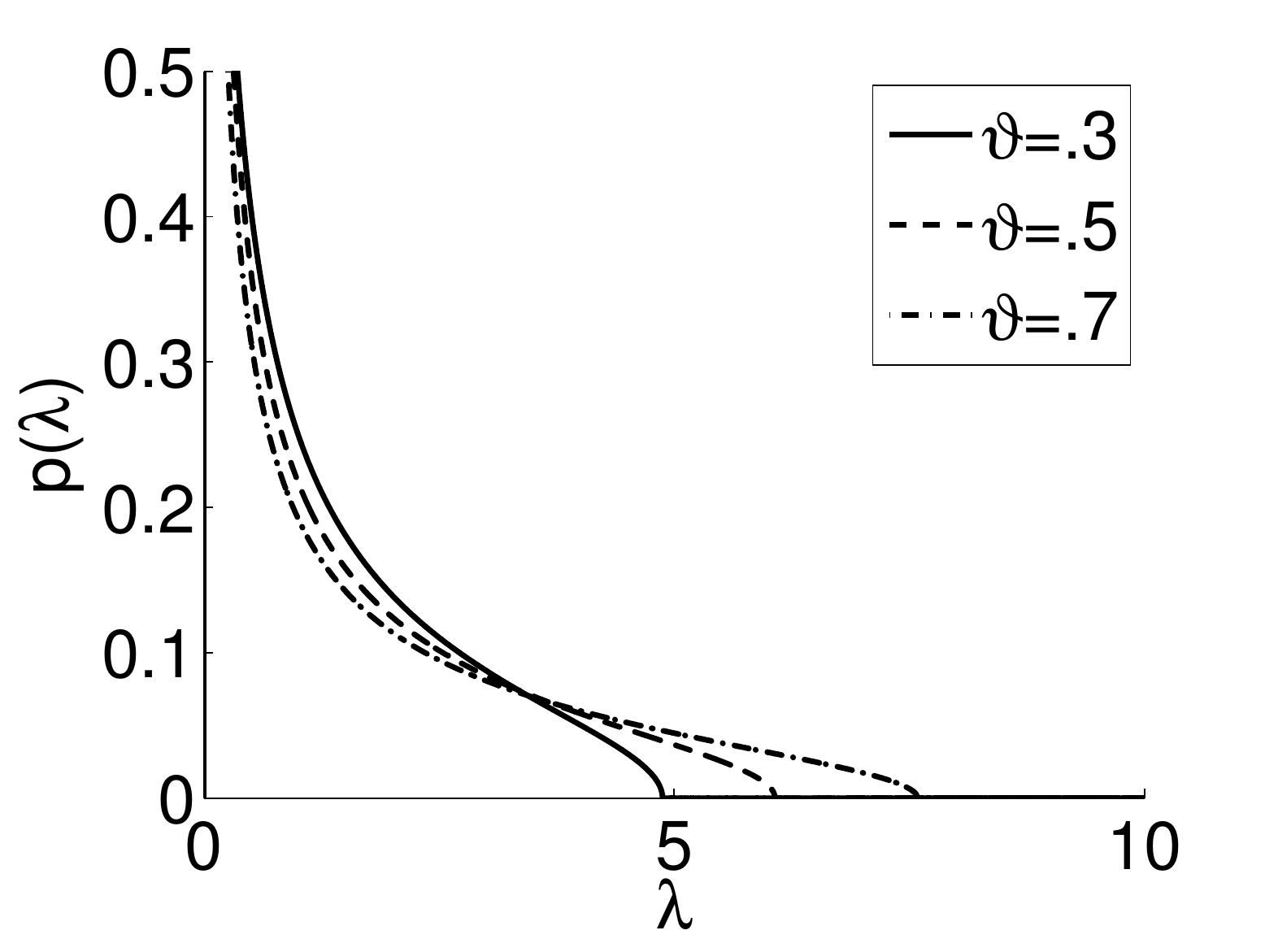}
}
\subfloat[$y=3$]{
\includegraphics[width=.31\linewidth]{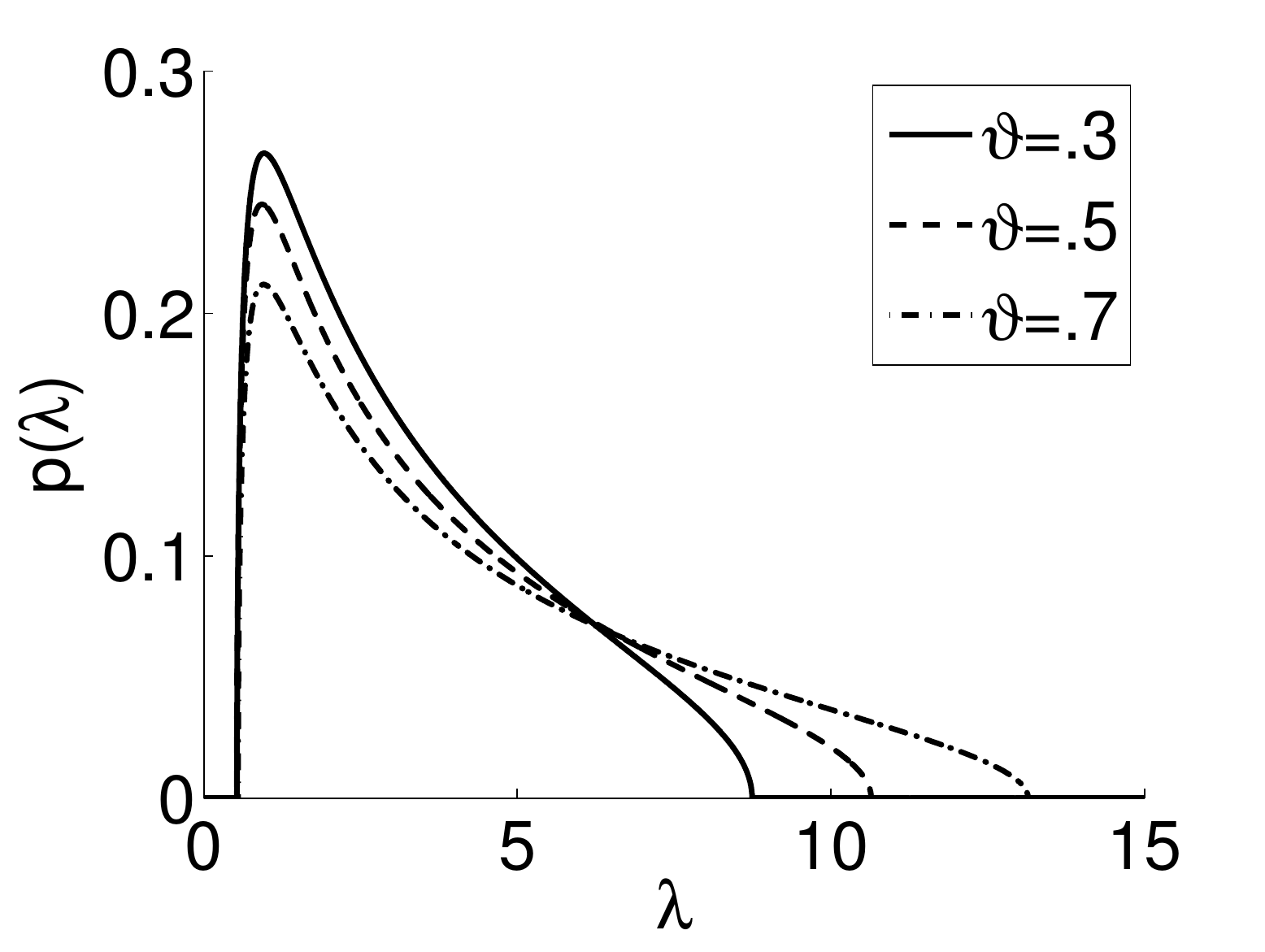}
}
\subfloat[$y=5$]{
\includegraphics[width=.31\linewidth]{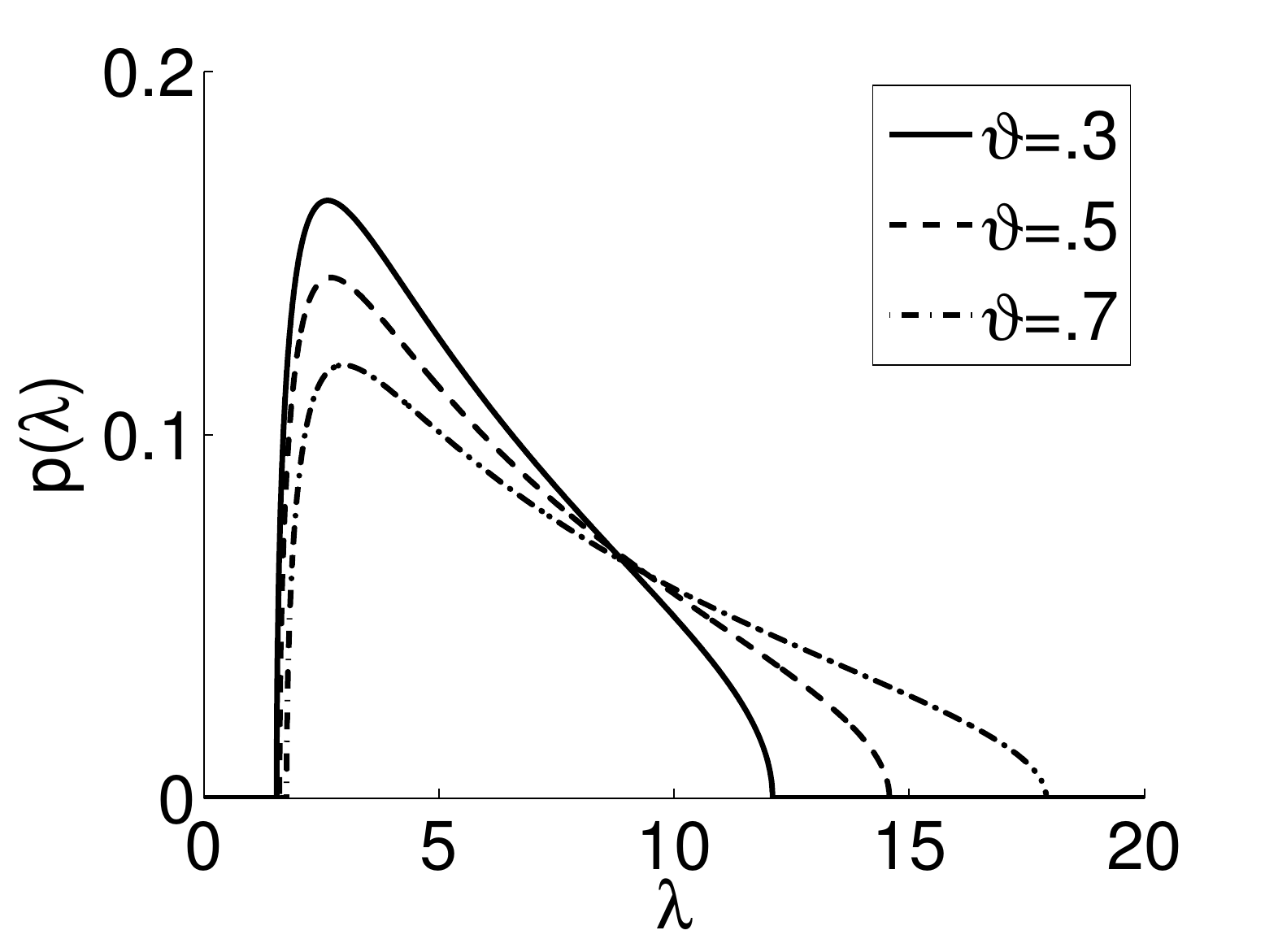}
}
\caption{Limiting spectral densities $\lambda\mapsto p(\lambda)$ of $p^{-1}\X\X^T$ for the MA(1) process $X_t=Z_t+\vartheta Z_{t-1}$ for different values of $\vartheta$ and $y=n/p$}
\label{fig-densitiesMA1}
\end{figure}

\begin{figure}
\subfloat[$y=1$]{
\includegraphics[width=.31\linewidth]{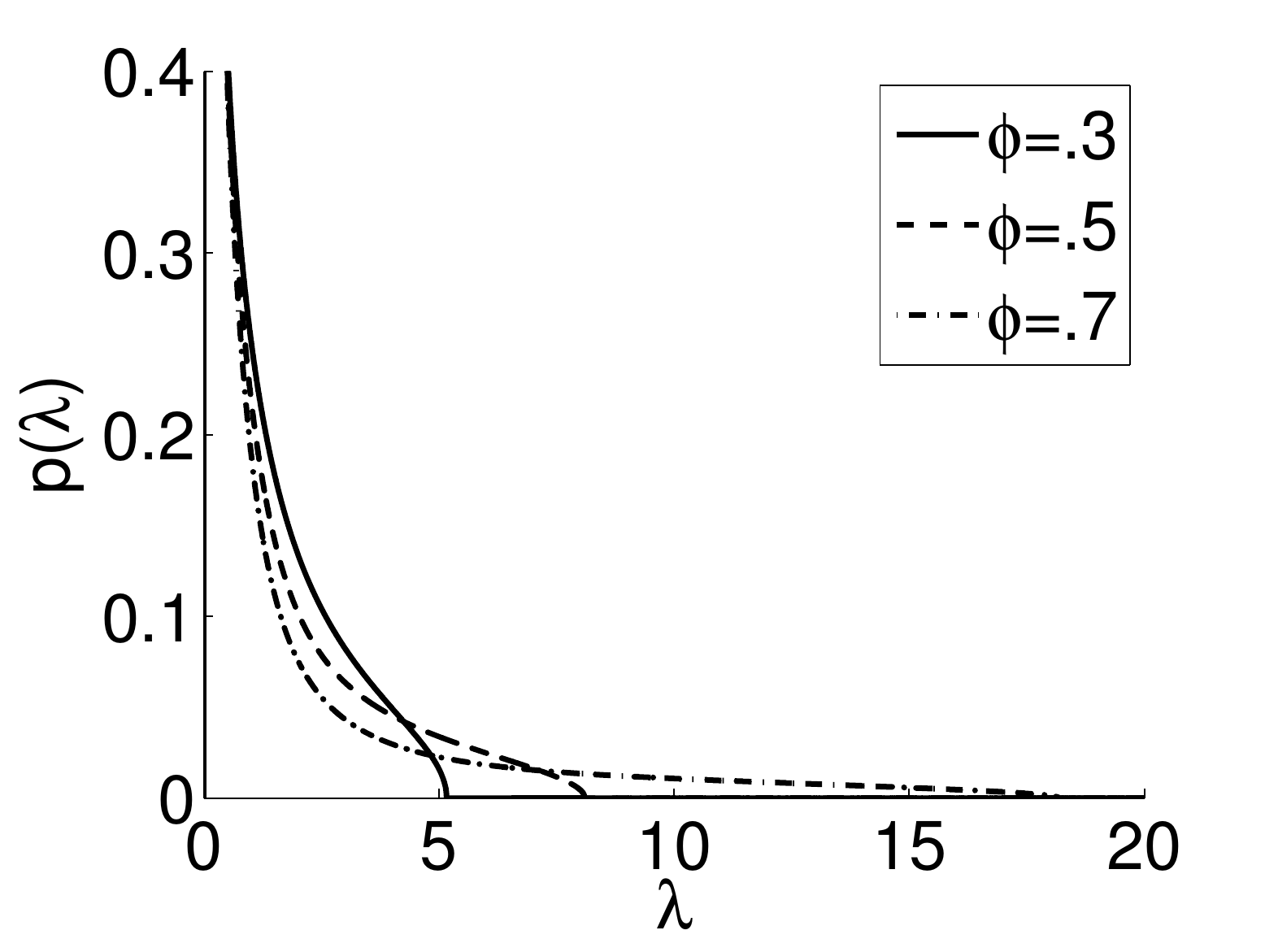}
}
\subfloat[$y=3$]{
\includegraphics[width=.31\linewidth]{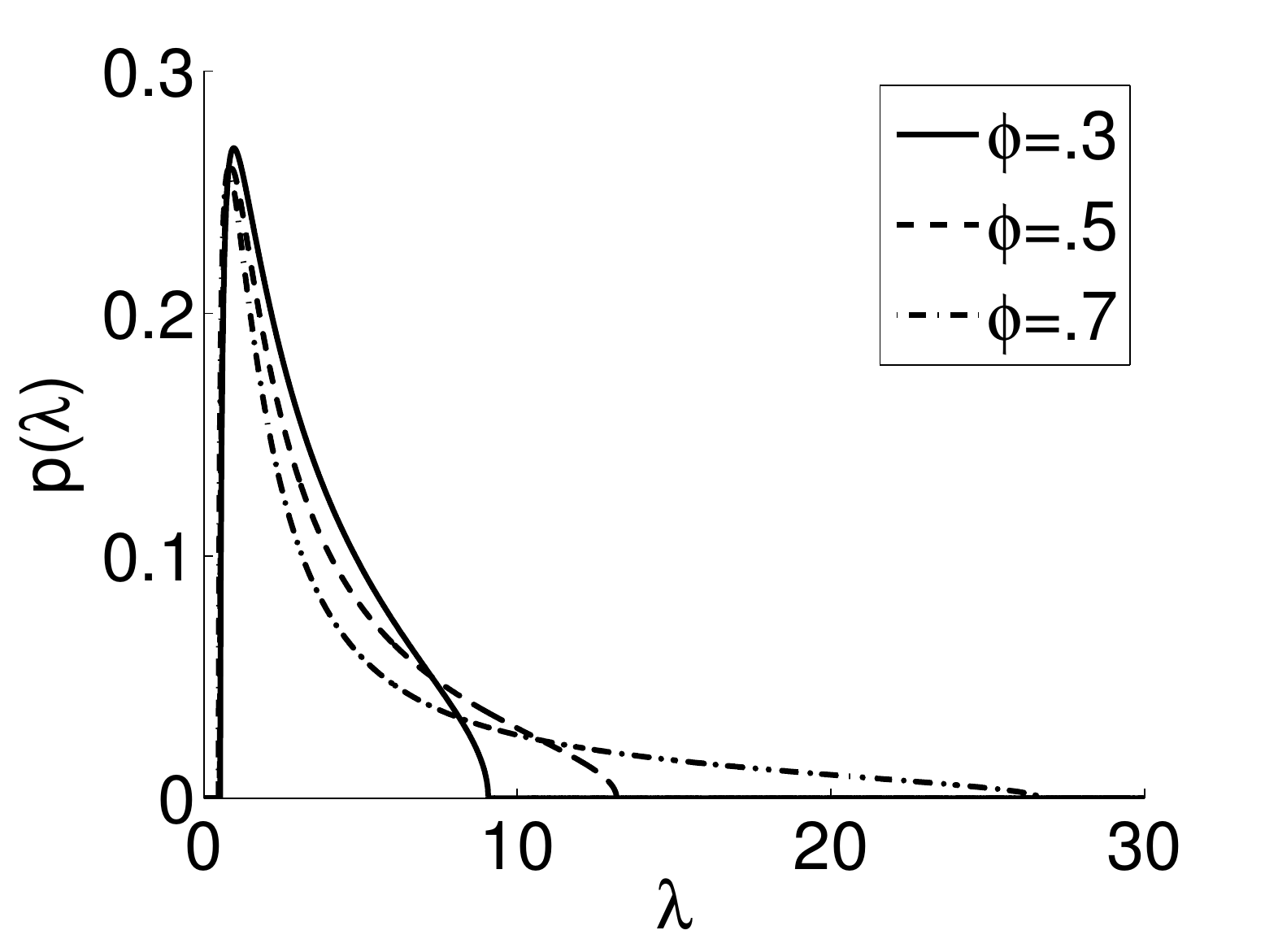}
}
\subfloat[$y=5$]{
\includegraphics[width=.31\linewidth]{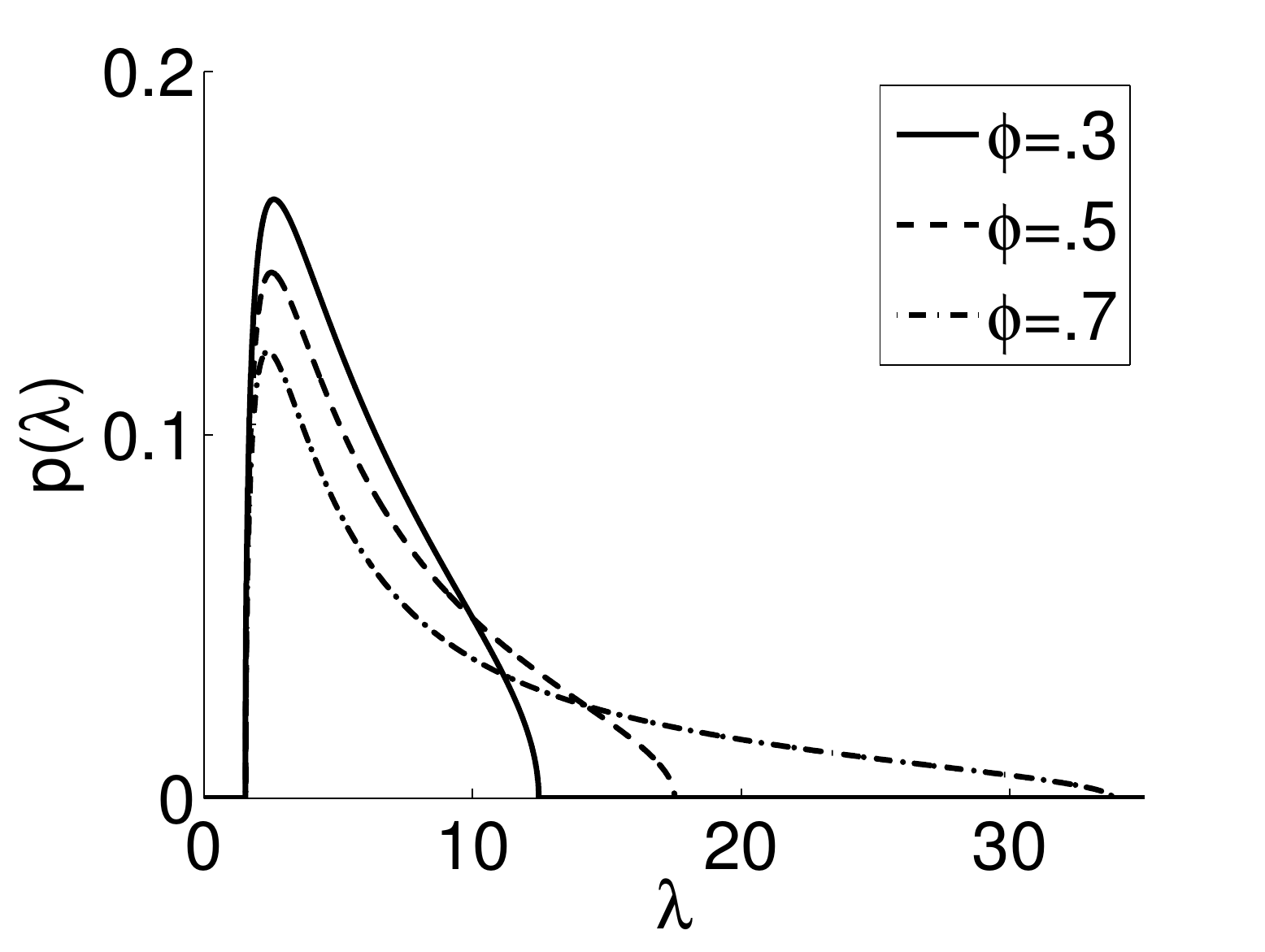}
}
\caption{Limiting spectral densities $\lambda\mapsto p(\lambda)$ of $p^{-1}\X\X^T$ for the AR(1) process $X_t=\varphi X_{t-1}+Z_t$ for different values of $\varphi$ and $y=n/p$}
\label{fig-densitiesAR1}
\end{figure}

\subsection{Autoregressive moving average processes}
Given polynomials $a:z\mapsto 1+a_1z+\ldots a_pz^p$ and $b:z\mapsto 1+b_1z+\ldots+b_qz^q$, an ARMA(p,q) process $X$ with autoregressive polynomial $a$ and moving average polynomial $b$ is defined as the stationary solution to the stochastic difference equation
\begin{equation*}
X_t+a_1X_{t-1}+\ldots+a_pX_{t-p} = Z_t+b_1Z_{t-1}+\ldots+b_qZ_{t-q},\quad t\in\Z.
\end{equation*}
If the zeros of $a$ lie outside the closed unit disk, it is well known that $X$ has an infinite\hyp{}order moving average representation $X_t=\sum_{j=0}^\infty c_j Z_{t-j}$, where $\{c_j\}$ are the coefficients in the power series expansion of $b(z)/a(z)$ around zero. It is also known (\citep[]{brockwell1991tst}) that there exist positive constants $\rho<1$ and $K$ such that $|c_j|\leq K\rho^j$, so that assumption \labelcref{maintheoremconditionsc} of \cref{maintheorem} is satisfied. While the autocovariance function of a general ARMA process does not in general have a simple closed form, its Fourier transform is given by
\begin{equation}
\label{eq-specdensARMA}
f(\omega)=\left|\frac{b\left(\ee^{\ii\omega}\right)}{a\left(\ee^{\ii\omega}\right)}\right|^2,\quad \omega\in[0,2\pi].
\end{equation}
Since $f$ is rational, assumptions \labelcref{maintheoremconditionsflevelsets,maintheoremconditionsdfzero} of \cref{maintheorem} are satisfied as well. In order to compute the LSD of $\Gamma$, it is necessary, by \cref{lemma-densityGamma}, to find the roots of a trigonometric polynomial of possibly high degree, which can be done numerically.

We now consider the special case of the ARMA(1,1) process $X_t=\varphi X_{t-1}+Z_t+\vartheta Z_{t-1}$, $|\varphi|<1$, for which one can obtain explicit results. By \cref{eq-specdensARMA}, the spectral density of X is given by
\begin{equation*}
f(\omega) = \frac{1+\vartheta^2+2\vartheta\cos\omega}{1+\varphi^2-2\varphi\cos\omega},\quad \omega\in[0,2\pi].
\end{equation*}
\Cref{eq-densityGamma} implies that the LSD of the autocovariance matrix $\Gamma$ has a density $g$, which is given by
\begin{align*}
g(\lambda)=&\frac{1}{2\pi}\sum_{\omega\in[0,2\pi]:f(\omega)=\lambda} \frac{1}{\left|f'(\omega)\right|}\\
  =&\frac{1}{\pi(\vartheta+\varphi\lambda)\sqrt{\left[(1+\vartheta)^2-\lambda(1-\varphi)^2\right]\left[\lambda(1+\varphi)^2-(1-\vartheta)^2\right]}}\1_{(\lambda_-,\lambda_+)}(\lambda),
\end{align*}
where
\begin{equation*}
\lambda_-=\min{(\lambda^-,\lambda^+)},\quad \lambda_+=\max{(\lambda^-,\lambda^+)},\quad \lambda^{\pm}=\frac{(1\pm\vartheta)^2}{(1\mp\varphi)^2}.
\end{equation*}
By \cref{maintheorem}, the Stieltjes transform $z\mapsto m_z$ of the limiting spectral distribution of $p^{-1}\X\X^T$ is the unique mapping $m:\C^+\to\C^+$ that satisfies the equation
\begin{align}
\label{eq-stieltjesARMA11}
\frac{1}{m_z} = & -z+y\int_{\lambda_-}^{\lambda_+}{\frac{\lambda g(\lambda)}{1+\lambda m_z}\dd\lambda}\notag\\
  = & -z+\frac{\vartheta y}{\vartheta m_z-\varphi}\\
    & -\frac{(\vartheta+\varphi)(1+\vartheta\varphi)y}{(\vartheta m_z-\varphi)\sqrt{\left[(1-\varphi)^2+m_z(1+\vartheta)^2\right]\left[(1+\varphi)^2+m_z(1-\vartheta)^2\right]}}.\notag
\end{align}
This is a quartic equation in $m_z\equiv m(z)$ which can be solved explicitly. An application of the Stieltjes inversion formula \labelcref{eq-Stieltjes} then yields the limiting spectral distribution of $p^{-1}\X\X^T$.

\begin{figure}
\subfloat[$y=1$]{
\includegraphics[width=.3\linewidth]{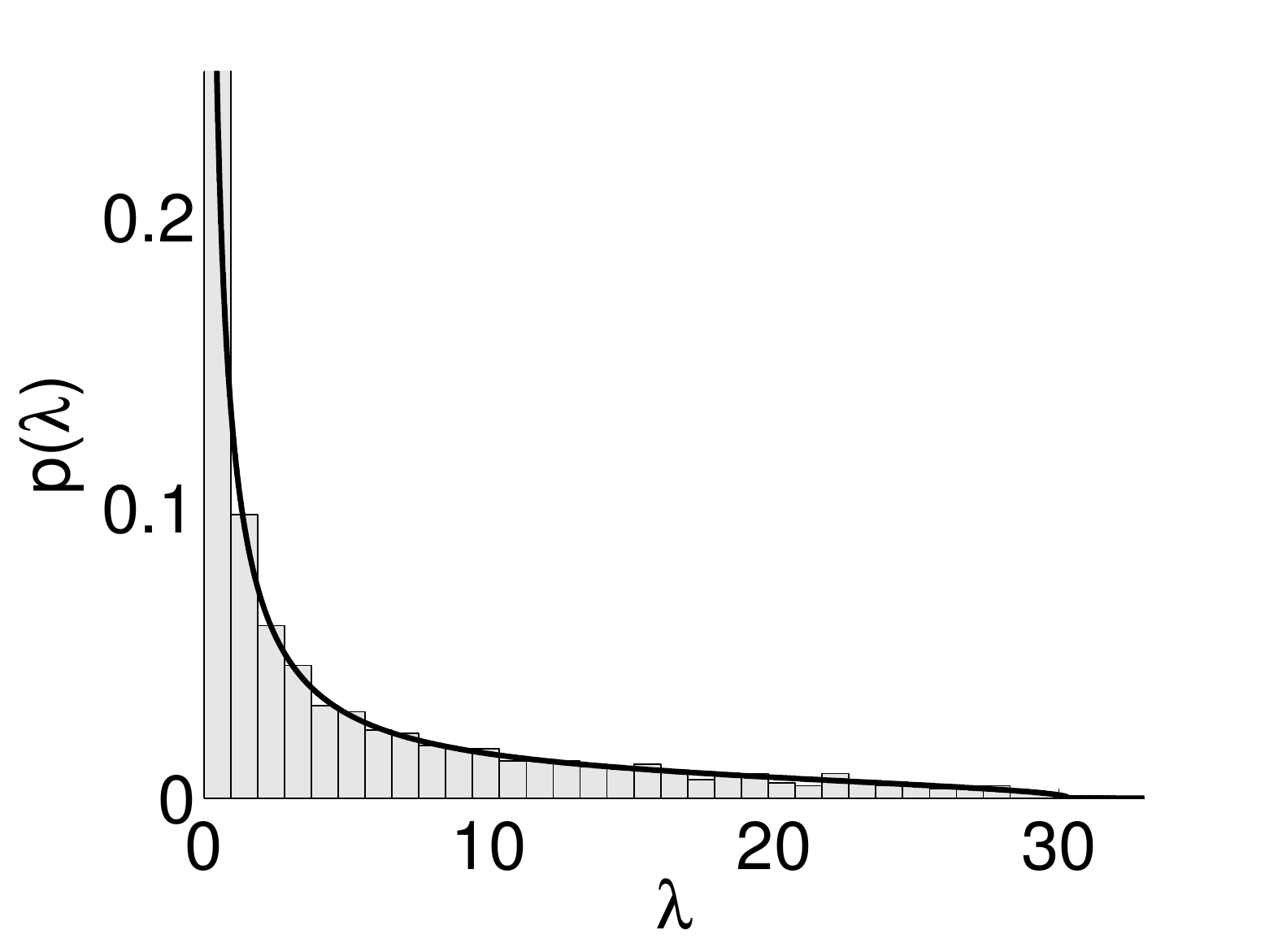}
}
\subfloat[$y=3$]{
\includegraphics[width=.3\linewidth]{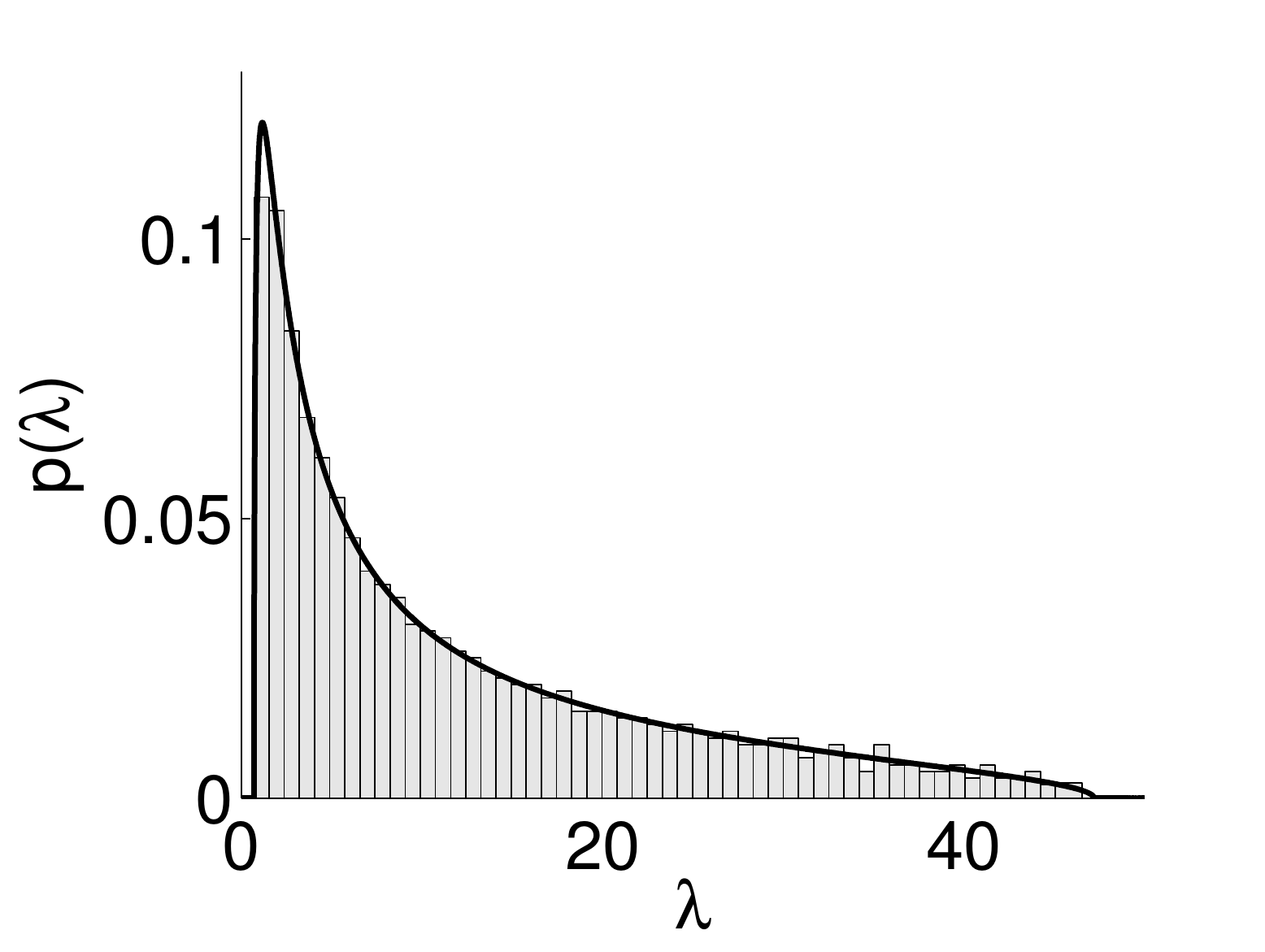}
}
\subfloat[$y=5$]{
\includegraphics[width=.3\linewidth]{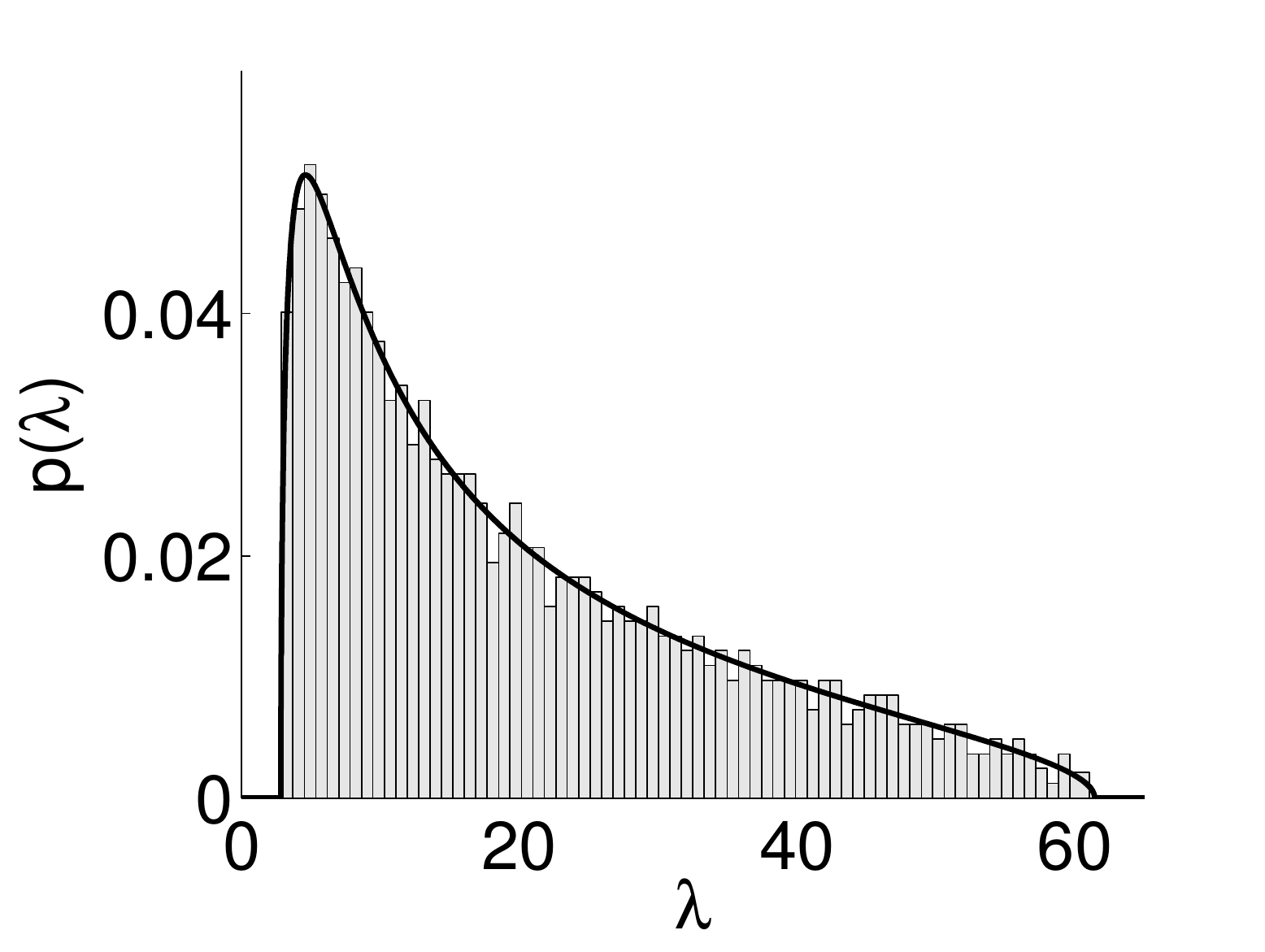}
}
\caption{Histograms of the eigenvalues and limiting spectral densities $\lambda\mapsto p(\lambda)$ of $p^{-1}\X\X^T$ for the ARMA(1,1) process $X_t=\frac{1}{2}X_{t-1}+Z_t+Z_{t-1}$ for different values of $y=n/p$, $p=1000$}
\label{fig-histoARMA11}
\end{figure}

If one sets $\varphi=0$, one obtains an MA(1) process; plots of the densities obtained in this case for different values of $\vartheta$ and $y$ are displayed in \cref{fig-densitiesMA1}. Similarly, the case $\vartheta=0$ corresponds to an AR(1) process; see \cref{fig-densitiesAR1} for a graphical representation of the densities one obtains for different values of $\varphi$ and $y$ in this case. For the special case $\varphi=1/2$, $\vartheta=1$, \cref{fig-histoARMA11} compares the histogram of the eigenvalues of $p^{-1}\X\X^T$ with the limiting spectral distribution obtained from \cref{maintheorem} for different values of $y$.

\Cref{eq-stieltjesARMA11} for the Stieltjes transform of the limiting spectral distribution of the sample covariance matrix of an ARMA(1,1) process should be compared to \citep[Eq. (2.10)]{bai2008}, where the analogous result is obtained for an autoregressive process of order one. They use the notation $c=\lim p/n$ and consider the spectral distribution of $n^{-1}\X\X^T$ instead of $p^{-1}\X\X^T$. If one observes that this difference in the normalization amounts to a linear transformation of the corresponding Stieltjes transform, one obtains their result as a special case of \cref{eq-stieltjesARMA11}.

\subsection{Fractionally integrated ARMA processes}
In many practical situations, data exhibit long\hyp{}range dependence, which can be modelled by long\hyp{}memory processes. Denote by $\BSO$ the backshift operator and define, for $d>-1$, the (fractional) difference operator by 
\begin{equation*}
\nabla^d=(1-\BSO)^d=\sum_{j=0}^\infty \prod_{k=1}^j \frac{k-1-d}{k} \BSO^j,\quad \BSO^jX_t=X_{t-j}.
\end{equation*}
A process $(X_t)_t$ is called a fractionally integrated ARMA(p,d,q) processes with $d\in(-1/2,1/2)$ and $p,q\in\NN$ if $(\nabla^d X_t)_t$ is an ARMA(p,q) process. These processes have a polynomially decaying autocorrelation function and therefore exhibit long\hyp{}range\hyp{}dependence, cf. \citep[Theorem 13.2.2]{brockwell1991tst} and \citep{granger1980,Hosking1981}. We assume that $d<0$, and that the zeros of the autoregressive polynomial $a$ of $(\nabla^d X_t)_t$ lie outside the closed unit disk. Then it follows that $X$ has an infinite\hyp{}order moving average representation $X_t=\sum_{j=0}^\infty c_j Z_{t-j}$, where the $(c_j)_j$ have, in contrast to our previous examples, not an exponential decay, but satisfy $K_1 (j+1)^{d-1} \leq c_j\leq K_2 (j+1)^{d-1}$, for some $K_1,K_2>0$. Therefore, if $d<0$, one can apply \cref{maintheorem} to obtain the LSD of the sample covariance matrix, using that the spectral density of $(X_t)_t$ is given by
\begin{equation*}
f(\omega)=\left|\frac{b\left(\ee^{\ii\omega}\right)}{a\left(\ee^{\ii\omega}\right)}\right|^2 \left|1-\ee^{-\ii\omega}\right|^{-2d},\quad \omega\in[0,2\pi].
\end{equation*}

\paragraph{Acknowledgements}
Both authors gratefully acknowledge financial support from Technische Universit\"at M\"unchen - Institute for Advanced Study funded by the German Excellence Initiative, and from the International Graduate School of Science and Engineering.

\end{document}